\documentclass[letterpaper, 10 pt, conference]{ieeeconf}

\IEEEoverridecommandlockouts           % This command is only
\let\labelindent\relax
\usepackage{url}
\usepackage{float}
\usepackage{enumitem}
\usepackage{graphicx}
\usepackage{amsthm,amssymb,amsfonts,thmtools,thm-restate}
\usepackage{grffile}
\usepackage{bigints}
\usepackage{tikz}
\usetikzlibrary{shapes.geometric, arrows, positioning, calc}
\usepackage{mathdots}
\usepackage{cite,cleveref}
\usepackage{caption,subcaption,multicol}
\usepackage{balance}
\usepackage{afterpage}
\newsavebox{\shortpagebox}

\makeatletter
\newcommand{\shortpage}[1]% #1= \twocolumn text to wrap into \onecolumn page
{\par
  \setbox\shortpagebox=\vbox{\strut #1\par}%
  \afterpage{\onecolumn
    \begin{multicols}{2}
    \unvbox\AP@partial
    \end{multicols}}%
  \unvbox\shortpagebox
\par}
\makeatother

\newcommand{\mcl}[1]{\mathcal{ #1}}
\newcommand{\mbf}[1]{\mathbf{ #1}}
\newcommand{\norm}[1]{\left\Vert #1\right\Vert}

\newcommand{\ip}[2]{\left\langle{#1},{#2}\right\rangle}

\newcommand{\bmat}[1]{\begin{bmatrix} #1\end{bmatrix}}

\newcommand{\mat}[1]{\begin{matrix}#1\end{matrix}}

\newcommand{\R}{\mathbb{R}}

\newcommand{\N}{\mathbb{N}}

\newtheorem{definition}{Definition}

\crefname{EqnBlock}{Figure}{Figures}

\newcommand{\threepi}[1]{\mcl{P}_{\{#1\}}}

\title{Computational stability analysis of PDEs with integral terms using the PIE framework}

\author{Sachin Shivakumar$^{1}$ \and Matthew M. Peet$^{1}$% <-this % stops a space
\thanks{$^{1}$ Sachin Shivakumar\{sshivak8@asu.edu\} and Matthew M. Peet\{mpeet@asu.edu\} are with School for Engineering of Matter, Transport and Energy, Arizona State University, Tempe, AZ, 85298 USA}%
}

\allowdisplaybreaks
\begin{document}
\maketitle

	%%%%%%%%%%%%%%%%%%%%%%%%%%%%%%%%%%%%%%%%%%%%%%%%%%%%%%%%%%%%%%%%%%%%%%%%%%%%%%%%
	\begin{abstract}
	The Partial Integral Equation (PIE) framework was developed to computationally analyze linear Partial Differential Equations (PDEs) where the PDE is first converted to a PIE and then the analysis problem is solved by solving operator-valued optimization problems. Previous works on the PIE framework focused on the analysis of PDEs with spatial derivatives up to $2^{nd}$-order. In this paper, we extend the class of PDEs by including integral terms and performing stability analysis using the PIE framework. More specifically, we show that PDEs with the integral terms where the integration is with respect to the spatial variable and the kernel of the integral operator is matrix-valued polynomials can be converted to PIEs if the boundary conditions satisfy certain criteria. The conversion is performed by using a change of variable where every PDE state is substituted in terms of its highest derivative and boundary values to obtain a new equation (a PIE) in a variable that does not have any continuity requirements. Later, we show that this change of variable can be represented using explicit maps from the parameters of the PDE to the parameters of the PIE and the stability test can be posed as an optimization problem involving these parameters. Lastly, we present numerical examples to demonstrate the simplicity and application of this method.
	\end{abstract}

	%%%%%%%%%%%%%%%%%%%%%%%%%%%%%%%%%%%%%%%%%%%%%%%%%%%%%%%%%%%%%%%%%%%%%%%%%%%%%%%%
	\section{Introduction}\label{sec:Introduction}

	Partial Differential Equations (PDEs) with integral terms (integrated with respect to the spatial variable) typically appear in the following two scenarios: a) when the inherent physics leads to integral terms while modeling (for example, population dynamics \cite{population_dynamics}, linear thermoelasticity \cite{day2013heat}, chemical reaction and transport \cite{chemical}, etc.); b) when the model involves point sensor measurements which, in practice, are weighted averages of the distributed state; and c) closed-loop dynamics of PDEs that are coupled with a boundary feedback controller or observer designed using the Backstepping methods \cite{krstic2008boundary}. Various mathematical analysis methods for such PDEs were presented in \cite{appell2000partial, gil2015stability} whereas \cite{yang2017stability} used computational tools like LMIs for analysis. However, these methods did not allow analysis (including stability analysis) of PDEs with integral terms at the boundary and typically can only be used for a fixed type of boundary conditions. 

	The goal of this paper is develop a computational method to analyze stability of PDEs with integral terms in both the dynamics and the boundary conditions. To describe the type of PDEs considered in this paper, we present the following two examples:

	\textit{Example 1:} Consider a model for population dynamics (called the McKendrick PDE; see \cite{population_dynamics}) given by 
	\begin{align*}
	\dot{x}(t,s) &= -\partial_s x(t,s) + f(s)x(t,s),\quad s\in[0,1],\\
	x(t,0) &= \int_0^1 h(s)x(t,s)ds,
	\end{align*}
	where $x$ is the density of population of age $s$ at any given time $t$. The boundary condition can be considered as an approximation of the number of newborns at time $t$ and is usually modeled as an integral boundary condition (a weighted average of the population density at different ages). 

	\textit{Example 2:} Consider the closed loop dynamics of a reaction-diffusion equation with a boundary sensed observer. The PDE dynamics along with the observer dynamics can be written as
	\begin{align*}
	&\dot{x}(t,s) = \lambda x(t,s) + \partial_s^2 x(t,s), \quad s\in [0,1],\\
	&\dot{\hat x}(t,s) = \lambda \hat x(t,s) + \partial_s^2 \hat x(t,s) + l(s) \left(\partial_s x(t,1) - \partial_s \hat x(t,1)\right),\\
	&\partial_s x(t,0) = 0, \quad x(t,1) = 0,\\
	&\partial_s \hat x(t,0) = 0, \quad \hat x(t,1) = 0.
	\end{align*}
	In \cite{krstic2008boundary}, the observer gains, $l$, were shown to be modified Bessel functions (denoted by $I_1$, a convergent but infinite series) and are given by 
	\begin{align*}
	l(s) =-\sqrt{\lambda} \frac{I_1\left(\sqrt{\lambda(1-s^2)}\right) }{\sqrt{1-s^2}}.
	\end{align*} A practical implementation of such an observer-based controller would involve an approximation of the kernels by a finite sum. While the original observer is provably stable, the stability of finite-dimensional approximation of observer is not guaranteed and must be verified --- a non-trivial step. Naturally, we want to simplify this post truncation verification step, possibly, by using a computational method. 

	Additionally, in practice, point measurements such as $y(t) = \partial_s x(t,1)$ (the sensor measurement that drives the observer dynamics) are not exact due to physical limitations of the sensing mechanism. A more accurate representation of such a sensor output is of the form
	\begin{align*}
	y(t) &= \int_0^1 w(s) \partial_s x(t,s) ds
	\end{align*}
	where $w$ is typically a Gaussian function with the center at $s=1$. Thus, a practical model of the observer dynamics can be written as
	\begin{align*}
	&\dot{\hat x}(t,s) = \lambda \hat x(t,s) + \partial_s^2 \hat x(t,s) \\
	&\quad +  \int_0^1l(s)w(\theta)\left(\partial_s x(t,\theta) - \partial_s \hat x(t,\theta)\right)d\theta.
	\end{align*}

	Note that Partial integrodifferential equations (PIDEs) such as Volterra integral equation and Fredholm integral equation are not considered since the variable of integration for such systems is the same as the variable in time-differential. To develop computational tests for the stability of PDEs with integral terms, We will use the Partial Integral Equation (PIE) framework, introduced in \cite{peet_2018CDC}, which is an approach to analysis/control PDEs by converting PDEs to PIEs and then solving analysis/control problems for the PIE using LMI based methods. This approach for analysis, although straightforward, cannot be directly used since a PIE form for PDEs with integral terms was not proposed in \cite{peet_2018CDC} or any of the subsequent works.

	The modified goal of this paper, then, is to extend the class of PDE models for which we have PDE to PIE conversion formulae to include the cases defined above. We approach the incorporation of a new class of PDE models in three steps: (a) we propose a parametric representation of the PDE model class and corresponding state-space; (b) We define an appropriate state-space to be used in the corresponding PIE model; (c) We find a unitary transformation from the PIE state-space to the state space of the PDE model -- proving equivalence of solutions and equivalence of stability properties. To summarize the main contribution of this paper, we provide explicit maps from the parameters of a PDE with integral terms to the parameters of the PIE where the PIE obtained will have stability properties identical to that of the PDE. Then we present the stability test for PDE with integral terms as an operator-valued optimization problem. 

	PIEs are explained in more detail in \cref{subsec:PIE}, however, here we briefly describe the PIE system and benefits of using the PIE framework to motivate our approach for stability analysis of PDEs with integral terms. A PIE: a) is defined by a single integro-differential equation; b) is parameterized by the $^*$ algebra of Partial Integral (PI) operators; and c) can be used to represent almost any well-posed PDE model. The general form of a PIE, defined by 3-PI operators, $\mcl T,\mcl A: L_2[a,b]\rightarrow L_2[a,b]$ is given by
	\begin{equation}\label{eq:PIE_intro}
	\mcl{T}\dot{\mbf{v}}(t) = \mcl{A}\mbf{v}(t)
	\end{equation}
	for all $t\ge 0$ where the state, $\mbf v(t) \in L_2[a, b]$ admits no continuity constraints or boundary conditions. An operator $\mcl P$ is said to be a 3-PI operator, denoted $\mcl P \in \Pi_3$ if there exist $R_0\in L_{\infty}$ and separable functions $R_1,R_2$, all matrix-valued, such that
	\begin{align*}
	\left(\mcl P\mbf u\right)(s) &= R_0(s)\mbf u(s) + \int_a^s R_1(s,\theta)\mbf u(\theta)d\theta \notag\\
	&\quad+ \int_s^b R_2(s,\theta)\mbf u(\theta) d\theta.
	\end{align*}

	The benefits of using a PIE for analysis, instead of a PDE, is because PIEs have a state-space form similar to ODEs and are parameterized by the PI operators that share many properties of matrices (specifically, being closed algebraically). Consequently, many numerical methods designed for control, analysis, and simulation of ODEs in state-space form can be extended to PIEs. For example, computational tests for analysis and control of PIEs were presented in \cite{peet_2020aut,shivakumar_2020CDC,robust_das2020_CDC}, etc. 

To summarize the organization of the paper, first, we define the algebra of PI operators in \Cref{subsec:PI} and define the class of PIE models, along with a definition of a solution in \Cref{subsec:PIE}. Then we introduce the class of PDE systems in \Cref{sec:PDE} -- including a definition of the solution of a PDE. Next, for a given PDE, we propose an admissibility condition that guarantees the existence of an equivalent PIE representation and present the formulae for constructing the PIE representation of the PDE (See \Cref{sec:pde2pie}). Finally, in \Cref{sec:numerical}, we consider several specific PDE models -- using the PIE representation and PIETOOLS to simulate the PDE and prove stability.	

\section{Notation}\label{sec:notation}
	We use the notation $0_{m\times n}$ to represent the zero matrix of dimension $m\times n$ and $0_n:=0_{n\times n}$. Similarly, $I_n$ is the identity matrix of dimension $n\times n$. We use $0$ and $I$ for the zero and identity matrix when the dimensions are clear from the context. $\R_+$ is the set of non-negative real numbers.

	$L_2^n[a,b]$ is the Hilbert space of $n$-dimensional vector-valued Lebesgue square-integrable functions defined on the interval $[a,b]$ and is equipped with the standard inner product. $L_{\infty}[a,b]^{m,n}$ is the Banach space of $m\times n$-dimensional essentially bounded measurable matrix-valued functions on $[a,b]$ equipped with the essential supremum singular value norm.

	Normal font $u$ or $u(t)$ typically implies that $u$ or $u(t)$ is a scalar or finite-dimensional vector (e.g. $u(t)\in \R^n$), whereas the bold font, $\mbf x$ or $\mbf x(t)$, typically implies that $\mbf x$ or $\mbf x(t)$ is a scalar or vector-valued function (e.g. $\mbf u(t)\in L_2^n[a,b]$). For a suitably differentiable function, $\mbf x$ of spatial variable $s$, we use $\partial_s^j \mbf x$ to denote the $j$-th order partial derivative $\frac{\partial^j\mbf x}{\partial s^j}$. For a suitably differentiable function of time and possibly space, we denote  $\dot{\mbf x}(t)=\frac{\partial }{\partial t}\mbf x(t)$. We use $W_{k}^n$ to denote the Sobolev spaces \begin{align}\label{eq:sobolev}
	W_{k}^n[a,b] := \{\mbf u\in L_2^n[a,b]\mid \partial^l_s \mbf u\in L_2^n[a,b] ~\text{for all} ~l\le k\}.
	\end{align}
	with inner product
	\[
	\ip{\mbf u}{\mbf v}_{W_{k}^n}=\sum\nolimits_{i=0}^k \ip{\partial_s^i\mbf u}{\partial_s^i \mbf v}_{L_2^n}.
	\]
	For brevity, we omit the domain $[a,b]$ and simply write $L_2^n$ or $W_{k}^n$ when clear from the context.

	% Finally, for normed spaced $A,B$, $\mcl L(A,B)$ denotes the space of bounded linear operators from $A$ to $B$ equipped with induced operator norm. $\mcl L(A) := \mcl L(A,A)$
	\subsection{PI Operators: A $^*$-algebra of bounded linear operators on $L_2$}\label{subsec:PI}
The PI algebras are parameterized classes of operators on $L_2^{n}$. The algebra of 3-PI operators which parameterize maps from $L_2^n \rightarrow L_2^n$.

\begin{definition}[Separable Function]\label{def:separable}
We say a function $R: [a,b]^{2} \rightarrow \R^{p \times q}$, is separable if there exist $r \in \N$ and functions $F\in L_{\infty}^{r\times p}[a,b]$  and $G\in L_{\infty}^{r\times q}[a,b]$  such that $R(s,\theta)=F(s)^T G(\theta)$.
\end{definition}

\begin{definition}[3-PI operators, $\Pi_3$]\label{def:3PI}
Given $R_0\in L_{\infty}^{p \times q}[a,b]$ and separable functions $R_1, R_2: [a,b]^{2} \rightarrow \R^{p \times q}$, we define the operator $\mcl{P}_{\{R_i\}}$ for $\mbf v \in L_2$ as
\begin{align}\label{eq:3pi}
\left(\mcl{P}_{\{R_i\}}\mbf v\right)(s)&:= R_0(s) \mbf v(s) +\int_{a}^s R_1(s,\theta)\mbf v(\theta)d \theta\\
&\quad +\int_s^bR_2(s,\theta)\mbf v(\theta)d \theta.
\end{align}
Furthermore, we say an operator, $\mcl P$, is 3-PI of dimension $p \times q$, denoted $\mcl P \in [\Pi_3]_{p,q}\subset \mcl L(L_2^q,L_2^p)$, if there exist functions $R_0\in L_{\infty}^{p\times q}$ and separable functions $R_1,R_2$ such that $\mcl P=\mcl{P}_{\{R_i\}}$.
\end{definition}

Additionally, we define $[\Pi_3]_{n,n}^+$ as the set of positive definite 3-PI operators where the positivity is defined with respect to $L_2$-inner product, i.e, $[\Pi_3]_{n,n}^+ := \{ \mcl P \in [\Pi_3]_{n,n}\mid \ip{x}{\mcl P x}_{L_2} > 0 ~\forall x\in L_2~\text{and}~x\ne 0\}$.

%\textbf{Note:} An important property of the PI algebras is that algebraic operations on the operator can be represented using algebraic operations on the parameters which represent that operator. 

\subsection{Partial Integral Equations}\label{subsec:PIE}
A Partial Integral Equation (PIE) is an extension of the state-space representation of ODEs (vector-valued first-order differential equations on $\R^n$) to spatially-distributed states on $L_2$. A PIE model is parameterized by 3-PI operators as
\begin{align}\label{eq:PIE_full}		
\mcl{T}\dot{\mbf{x}}_f(t)&=\mcl{A}\mbf x_f(t), \mbf{x}_f(0) = \mbf{x}_f^0\in L_2^{m,n}[a,b],
\end{align}
where $\mcl T\in [\Pi_3]_{n,n}$ and $\mcl A\in [\Pi_3]_{n,n}$ are 3-PI operators. 

Unlike PDE models, a PIE does not allow for spatial derivatives -- only a first-order derivative with respect to time. In particular, the state of the PIE model, $\mbf x_f \in L_2[a,b]$ is not differentiable and consequently, no boundary conditions are possible in the PIE framework. Finally, given a PIE, we require differentiability of the solution with respect to the $\mcl{T}$-norm which is defined as
\begin{align}\label{eq:Tnorm}
\norm{\mbf{x}_f}_{\mcl{T}} := \norm{\mcl{T}\mbf{x}_f}_{L_2}, \qquad \text{for}~\mbf{x}_f\in L_2.
\end{align}

	\begin{definition}[Definition of solution of a PIE]\label{def:piesolution}
	For given inputs $\mbf x_f^0\in L_2^{n}$, we say that $\{\mbf x_f\}$ satisfies the PIE defined by $\{\mcl{T},$ $\mcl{A}\}$ with initial condition $\mbf x_f^0$ if
	\begin{enumerate}
		\item $\mbf x_f(t)\in L_2^{n}[a,b]$ for all $t\ge 0$
		\item $\mbf x_f$ is Frech\'et differentiable with respect to the $\mcl{T}$-norm almost everywhere on $\R_+$
		\item $\mbf x_f(0) = \mbf x_f^0$
		\item Eq.~\eqref{eq:PIE_full} is satisfied for almost all $t\in\R_+$.
	\end{enumerate}
	\end{definition}

\section{A parametric form of PDEs with integral terms}\label{sec:PDE}

We divide the parameters of a PDE into three groups depending based on the type of constraint they appear in -- namely the continuity constraints, the in-domain dynamics, and the boundary conditions. The continuity constraints specify the existence of partial derivatives and boundary values for each state as required by the in-domain dynamics and boundary conditions. The boundary conditions are represented as a real-valued algebraic constraint that maps the distributed state to a vector of boundary values which are then used to constrain the in-domain dynamics. The in-domain dynamics (or generating equation) specify the time derivative of the state at every point in the interior of the domain and allow for both integral and spatial derivative operators. In the following subsections, we highlight the parameters required to unambiguously define the above three types of constraint in a PDE.

	\subsubsection{The continuity constraint} Given a PDE state, $\mbf x(t,\cdot)$, the continuity constraint can be uniquely defined by the parameter $n=\{n_0,n_1,n_2\}$, which partitions and orders the PDE states $\mbf x$ by increasing differentiability as follows.
	\[
	\mbf x(t,\cdot)=\bmat{\mbf x_0(t,\cdot)\\ \mbf x_1(t,\cdot)\\\mbf x_2(t,\cdot)} \in \bmat{W_0^{n_0}\\W_1^{n_1} \\ W_2^{n_2}} =: W^n.
	\]
	Given such an $n\in \N^3$, we can identify all well-defined partial derivatives of $\mbf x$.

	\noindent \textbf{Notation:} For convenience, we define the vector of all continuous partial derivatives of the PDE state $\mbf x$ as permitted by the continuity constraint as $\mbf x_c$, the vector of all partial derivatives as $\mbf x_D$ and the list of all possible boundary values of $\mbf x$ as $x_b$, i.e.,
	\begin{flalign}\label{eq:odepde_bc}
		%%%%%%%%%%%%%%%%%%%%%%%%%%%%%%%%%%%
		\mat{\mbf x_c(t,\cdot)=\bmat{\mbf x_1(t,\cdot)\\\mbf x_2(t,\cdot)\\\partial_s \mbf x_2(t,\cdot)}\\\\\quad x_b(t)=\bmat{\mbf x_c(t,a)\\\mbf x_c(t,b)}},\quad \mbf x_D(t,\cdot) = \bmat{\mbf x_0(t,\cdot)\\\mbf x_1(t,\cdot)\\\mbf x_2(t,\cdot)\\\partial_s\mbf x_1(t,\cdot)\\\partial_s\mbf x_2(t,\cdot)\\\partial_s^2 \mbf x_2(t,\cdot)}.
		\end{flalign}	
Additionally, we define some notation given the continuity constraint parameter $n=\{n_0,n_1,n_2\}$ that will repeatedly appear in the subsequent sections. Let $n_{\mbf x}:=\sum_{i=0}^2 n_i$ be the number of states in $\mbf x$, $n_{S_i}:=\sum_{j=i}^2 n_{j}$ be the number of $i$-times differentiable states, and $n_{S}=\sum_{i=1}^2 n_{S_i}$ be the total number of possible partial derivatives of $\mbf x$.
		\subsubsection{Boundary Conditions}
		Given an $n=\{n_0,n_1,n_2\}$, we now parameterize a generalized class of boundary conditions consisting of a combination of boundary values and integrals of the PDE state. Specifically, the boundary conditions are parameterized by square integrable functions $B_{I,i}(s) \in \R^{n_{BC} \times n_{Si}}$ for $i\ge 0$ and $B\in \R^{n_{BC}\times 2n_S}$ as
		\begin{flalign}\label{eq:odepde_b}
		0 &= \int_{a}^{b} B_{I}(s)\mbf x_D(t,s)ds -Bx_b(t)
		\end{flalign}
		where $n_{BC}$ is the number of specified boundary conditions. For reasons of well-posedness, as discussed in \Cref{sec:pde2pie}, we typically require $n_{BC}=n_{S}$. 

		These boundary conditions and the continuity constraints collectively define the domain of the infinitesimal generator -- which specifies a set of acceptable solutions $\mbf x(t)\in X$ for the PDE -- as 
		\begin{align}\label{eq:odepde_general_domain}
		&X = \notag\\
		&\left\lbrace \mbf{x}\in W^{n}[a,b]: B\bmat{\mbf x_c(t,a)\\\mbf x_c(t,b)}=\int\limits_{a}^{b} B_{I}(s)\mbf{x}_D(t,s)ds\right\rbrace.
		\end{align}

		\noindent\textbf{Notation:} For convenience, we collect all the parameters which define the boundary-valued constraint in Eq.~\eqref{eq:odepde_b} and introduce the shorthand notation $\mbf G_{\mathrm{b}}$ which represents the labelled tuple of system parameters as
		\begin{align}
		&\mbf G_{\mathrm{b}} = \left\{B,~B_{I}\right\}.\label{eq:BC-parms}
		\end{align}
		When this shorthand notation is used to denote a given set of system parameters, it is presumed that all parameters have appropriate dimensions.

		\subsubsection{Dynamics of the PDE}
		Finally, we may now define the dynamics of the PDE which is parameterized by the functions $A_{0}(s),$ $A_{1}(s,\theta),$ $A_{2}(s,\theta)$ $\in \R^{n_{\mbf x} \times n_{S}}$, and $B_{xb}(s)$ $\in\R^{n_{\mbf x}\times 2n_S}$ as
		\begin{flalign}\label{eq:general_pde_}
		\dot{\mbf{x}}(t,s) &= A_0(s)\mbf x_D(t,s) +\int\limits_{a}^{s}A_1(s,\theta)\mbf x_D(t,\theta)d\theta\notag\\
		&\quad+\int\limits_{s}^{b}A_2(s,\theta)\mbf x_D(t,\theta)d\theta,
		\end{flalign}
		with the constraint $\mbf x(t)\in X$. In this representation, the kernels $A_{1},A_{2}$ represent integral operators.

		\noindent\textbf{Notation:} For convenience, we collect all the parameters from the generating equation (Eq.~\eqref{eq:general_pde_}) and introduce the shorthand notation $\mbf G_{\mathrm{p}}$ which represents the labelled tuple of system parameters as
		\begin{align}
	%&\mbf G_{\mathrm{o}} = \left\{A,~B_{xi},~C_{zx},~D_{zi},~C_{yx},~D_{yi},~C_v,~D_{vw},~D_{vu}\mid i\in \{x,w,u\}\}\right\},\label{eq:ode-params}\\
	&\mbf G_{\mathrm{p}} = \left\{A_{0},~A_{1},~A_{2}\right\}.\label{eq:pde-parms}
	\end{align}
	When this shorthand notation is used to denote a given set of system parameters, it is presumed that all parameters have appropriate dimensions. Now, we define the notion of a solution to the above described PDE .

	\begin{definition}[Definition of solution of a PDE]\label{defn:PDE}
	For given $\mbf x^0 \in X$, we say that $\{\mbf x\}$ satisfies the PDE defined by $\{n,\mbf G_{\mathrm b}, \mbf G_{\mathrm{p}}\}$ (defined in Eqs.~\eqref{eq:BC-parms} and \eqref{eq:pde-parms}) with initial condition $\mbf{x}^0$ if
	\begin{itemize}
		\item $\mbf x(t) \in X$ for all $t\ge 0$
		\item $\mbf x$ is Frech\'et differentiable with respect to the $L_2$-norm almost everywhere on $\R_+$
		\item $\mbf x(0) = \mbf x^0$
		\item Eq.~\eqref{eq:general_pde_} is satisfied for almost all $t\ge 0$.
	\end{itemize}
	\end{definition}

\section{Representing a PDE as a PIE}\label{sec:pde2pie}
	In \Cref{sec:PDE}, we presented the PDE form Eq.~\eqref{eq:general_pde_}. Recall that our goal is to find a PIE form of such PDEs so that the computational tools developed for PIEs can be used in analysis/control. Specifically, we will express PDEs introduced in previous section as PIEs of the form
	\begin{align}\label{eq:PIE_}			
	\mcl{T}\dot{\mbf{x}}_f(t)&=\mcl{A}\mbf x_f(t),\qquad \mbf{x}_f(0) = \hat{\mbf{x}}_f^0\in L_2^{n}.
	\end{align}
	However, first, we have to verify that such a PIE form does indeed exist, which can be verified by testing some additional constraints on $\mbf G_{\mathrm b}$. Specifically, we can verify that for any well-posed PDE of the form given in \Cref{sec:PDE} with the initial condition $\mbf x^0\in X$, there exists a corresponding PIE with corresponding initial condition $\mbf x_f^0 \in L_{2}^n$ whose solution can be used to construct a solution to the PDE --- i.e., there exists a invertible map from any initial condition and corresponding solution of the PDE to a corresponding initial condition and solution of the corresponding PIE.
	In the following subsection, we propose a sufficient condition for the existence of such an invertible map. 
	
	\subsection{Admissibility of the Boundary Conditions}\label{subsec:assumptions}
	The following definition of admissibility imposes a notion of well-posedness on $X$, the domain of the PDE defined by the continuity constraints and the boundary conditions which when satisfied guarantees the existence of a PIE form for the PDE.
	\begin{definition}[Admissible Boundary Conditions]
	Given an $n$ and a parameter set, $\mbf G_{\mathrm b}$, we say the pair $\{n,\mbf G_{\mathrm b}\}$ is \textbf{admissible} if $B_T$ is invertible where
	\begin{align}\label{as:invertible}
	&B_T := B\bmat{T(0)\\T(b-a)}-\int_{a}^{b}B_{I}(s)U_2 T(s-a)ds,
	\end{align}
	and $B_I$, $T$, and $U_1$ are as defined in \Cref{fig:Gb_definitions}.
	\end{definition}
	Since $B_T$ is invertible only when it is a square matrix, naturally, we require $n_{BC}=n_S$, i.e., when we have $n_S$ differentiable states we need $n_S$ boundary conditions for a well-posed solution. Note that the test for invertibility of $B_T$ depends only on the boundary condition parameters and not the dynamics or the initial condition.
	Hence, we refer to this test as ``well-posedness of the boundary conditions'', which is not necessarily the same as ``well-posedness of a PDE''. 
	
	\subsection{PIE representation of a PDE: $\mcl T$ and $\mcl A$}

	Assuming that we have a PDE with admissible $\{n,\mbf G_{\mathrm b}\}$, in this subsection, we propose a PIE form for any PDE of the form Eq.~\eqref{eq:general_pde_} such that the solutions of the PDE and the corresponding PIE are equivalent, in the sense, the existence of the PDE solution guarantees the existence of a PIE solution and vice versa. The conversion of PDE to PIE and the equivalence of the solutions are summarized below.

\begin{restatable}{theorem}{thmEquivalence}\label{thm:equivalence}
Given a set of PDE parameters $\{n,$ $\mbf G_{\mathrm b},$ $\mbf G_{\mathrm{p}}\}$ as defined in \Cref{eq:pde-parms,eq:BC-parms} with $\{n,\mbf G_{\mathrm b}\}$ admissible, let the PI operators $\{\mcl{T},$ $\mcl A\}$ be as defined in \cref{fig:Gb_definitions}. Then, for any $\mbf{x}^0$ $\in$ $X$ ($X$ is as defined in \Cref{eq:odepde_general_domain}), $\{\mbf{x}\}$ satisfies the PDE defined by $\{n,\mbf G_{\mathrm b}, \mbf G_{\mathrm{p}}\}$ with initial condition $\mbf{x}^0$ if and only if $\{\mcl D\mbf{x}\}$ satisfies the PIE defined by $\{\mcl{T},$ $\mcl A\}$ with initial condition $\mcl D\mbf{x}^0 \in L_2^{n_{\mbf x}}$ where $\mcl D\mbf x = \text{col}(\partial_s^0\mbf x_0,\partial_s\mbf x_1,\partial_s^2\mbf x_2)$.
\end{restatable}
\begin{proof}
The proof is similar to the proof of Theorem 4.4 and 4.5 in \cite{shivakumar_SICON2022} with $n=\{n_0,n_1,n_2\}$, $\hat{\mcl T} = \mcl T$, $\hat{\mcl A} = \mcl A$ and other unused parameters set to empty sets or zeros.
\end{proof}

The above result provides a map from the PDE solution, $\mbf x$, to the solution of the corresponding PIE, $\mcl D\mbf x$. However, there also exists an inverse map which is given by $\mcl T$ that maps the solution of the PIE back to the PDE solution. The invertible relation between the solutions can be summarized as follows. 

% \begin{theorem}
\begin{restatable}{theorem}{Tmap}\label{thm:T_map}
Given an $n$, and $\mbf G_{\mathrm b}$ with $\{n,\mbf G_{\mathrm b}\}$ admissible,  let $\mcl{T}$ be as defined in \cref{fig:Gb_definitions}, $X$ as defined in Eq.~\eqref{eq:odepde_general_domain} and $\mcl D$ $:=$diag$(\partial_s^0 I_{n_0},$ $\partial_s I_{n_1},$ $\partial_s^2 I_{n_2})$. Then we have the following.
\begin{enumerate}[label=(\alph*)]
\item If $\mbf{x} \in X$, then $\mcl D \mbf x \in L_2^{n_{\mbf x}}$ and $\mbf{x} = \mcl{T}\mcl D \mbf x$.
\item For any $\hat{\mbf x} \in L_2^{n_{\mbf x}}$, $\hat{\mcl{T}}  \hat{\mbf{x}}\in X$ and $\hat{\mbf{x}} = \mcl D \mcl{T}  \hat{\mbf{x}}$.
\end{enumerate}
\end{restatable}
% \end{theorem}	

\begin{proof}
The proof is similar to the proof of Theorem 4.3 in \cite{shivakumar_SICON2022} with $n=\{n_0,n_1,n_2\}$, $\hat{\mcl T} = \mcl T$ and other unused parameters set to empty sets or zeros.
\end{proof}
This bijective mapping ensures the existence of both solutions when one of them exists. Given the PDE parameters $\{n, \mbf G_{\mathrm b}, \mbf G_{\mathrm p}\}$ we now have a set of formulae to find the corresponding PIE parameters $\{\mcl T,\mcl A\}$. Furthermore, we know that a solution to the PDE equation provides a unique solution to the PIE and vice versa.

% \subsection{PIE representation of a PDE: $\{\mcl T,$ $\mcl A\}$}\label{subsec:pde_equivalence}
% Now that we have found an invertible state transformation of the PDE ($\mbf x \rightarrow \hat{\mbf x}_f={\mcl D} \mbf x$), we construct the corresponding PIE representation by replacing all $\mbf x$ terms in the PDE with $\mbf x=\mcl T (\mcl D\mbf x)=\mcl T \mbf x_f$. This substitution results in a PIE of the form Eq.~\eqref{eq:PIE_} where the relevant PI operators are given in \cref{fig:Gb_definitions}.
% Note that, in the formulae, we have used the notation $\mbf P_c$ which stands for the composition operation of two PI operators. 

\section{Equivalence of representations: PDE and PIE}\label{sec:equivalence}
	In this section, we show that the solutions to the two representations (PDE and PIE) have equivalent stability properties and present a solvable optimization problem to prove stability. However, we have to first define a notion of stability for the two representations.
	%%%%%%%%%%%%%%%%%%%%%%%%%%%%%%%%%%%%%%%%%%%%%%%%%%%%%%%%%%%%%%%%%%%%%%%%%%%%%%%%%%%%
	%\newpage
	\subsection{Definitions of stability}\label{sec:stab_pde}
	We define the notion of stability of an PDE based on the stability of $\mbf{x}$ that satisfies the PDE for given initial conditions where the stability is defined with respect to the standard Sobolev norm $H$ which is defined as
	\begin{align*}
	\norm{\mbf x}_{H} = \sum_{i=0}^2 \norm{\mbf x_i}_{H_i},\quad \norm{\mbf x_i}_{H_i}: = \sum_{j=0}^{i}\norm{\partial_s^j \mbf x_i}_{L_2}.
	\end{align*}
	
	\begin{definition}[Exponential Stability of a PDE]
	We say a PDE defined by $\{n, \mbf G_{\mathrm b}, \mbf G_{\mathrm p}\}$ is exponentially stable, if there exists constants $M$, $\alpha>0$ such that for any $\mbf{x}^0\in X$, if $\mbf{x}$ satisfies the PDE defined by $\{n, \mbf G_{\mathrm b}, \mbf G_{\mathrm p}\}$ with initial condition $\mbf x^0$ then
	\begin{align*}
	\norm{\mbf{x}(t)}_{H}\le M\norm{\mbf{x}^0}_He^{-\alpha t} \qquad \text{for all} ~t\ge 0.
	\end{align*}
	\end{definition}
	
	Similar to the stability of a PDE, we can define the stability of a PIE system based on stability $\mbf x_f$ that satisfies the PIE for some initial condition and zero inputs. Unlike, the stability of PDE, the stability of a PIE is defined with respect to the $L_2$-norm since solutions of PIE need not have spatial continuity.
	
	\begin{definition}[Exponential Stability of a PIE]
	We say a PIE defined by $\{\mcl T,\mcl A\}$ is exponentially stable, if there exists constants $M$, $\alpha>0$ such that for any $\mbf{x}_f^0\in L_2^{n}$, if $\mbf{x}_f$ satisfies the PIE defined by $\{\mcl T,\mcl A\}$ with initial condition $\mbf x^0_f$, then
	\begin{align*}
	\norm{\mbf{x}_f(t)}_{ L_2}\le M\norm{\mbf x_f^0}_{ L_2}e^{-\alpha t} \qquad \text{for all} ~t\ge 0.
	\end{align*}
	\end{definition}
	
	Based on these definitions, the goal now is to show that for any $\mbf x$ that satisfies the PDE and $\mbf x_f$ that satisfies the PIE, if $\mbf x = \mcl T\mbf x_f$ with $\mcl T$ invertible, then exponential decay of $\mbf x$ implies exponential decay of $\mbf x_f$ and vice versa. The main hurdle in proving this is the fact that these two notions of stability are defined using different norms ($\norm{\cdot}_H$ and $\norm{\cdot}_{L_2}$). For any $\mbf x$ and $\mbf x_f$ such that $\mbf x = \mcl T\mbf x_f$, we know that $\norm{\mbf x}_{H}\le c$ implies $\norm{\mcl T\mbf x_f}_{L_2}<c$, but the converse is typically not true. Thus, additional steps are needed to prove the converse implication which is accomplished using a new norm on the space $X$, denoted by $\norm{\cdot}_X$, to show that:
	\begin{enumerate}
		\item $\mcl T$ is a norm-preserving bijection from $L_2$ to $X$ (when equipped with $\norm{\cdot}_X$).
		\begin{itemize}
			\item[-] Consequently, $X$ is closed under $\norm{\cdot}_X$ (and $\mcl T$ is unitary)
		\end{itemize}
		\item $\norm{\cdot}_H$ is equivalent to $\norm{\cdot}_X$ (and thus also equivalent to $\norm{\cdot}_{L_2}$) on the subspace $X$
		\item Finally, the PDE $\{n,\mbf G_{\mathrm p},\mbf G_{\mathrm b}\}$ is exponentially stable if and only if the corresponding PIE defined by $\{\mcl T,\mcl A\}$ is exponentially stable (PIE parameters are related to PDE parameters as shown in \cref{fig:Gb_definitions}).
	\end{enumerate}

	\subsection{The map $\mcl T$ are unitary}\label{subsec:unitary}
	First, we would like to show that $X$ is complete with respect to the norm $\norm{\cdot}_{X}$. Previously, in \Cref{thm:T_map}, we showed that $\mcl T$ is invertible. Therefore, we just need to show that $\mcl T$ preserves the inner product. However, we first define the new $X$-inner product as
	\begin{align}\label{eq:ip_xv}
	\ip{\mbf{x}}{\mbf{y}}_{X}:= \sum_{i=0}^2 \ip{\partial_s^i\mbf{x}_i}{\partial_s^i\mbf{y}_i}_{L_2}= \ip{\mcl D\mbf x}{\mcl D \mbf y}_{L_2}.
	\end{align}	
	
	\begin{restatable}{theorem}{thmUnitary}\label{thm:unitary_T}
	Suppose $\{n,\mbf G_{\mathrm b}\}$ is admissible, $\mcl T$ is as defined in \Cref{fig:Gb_definitions} and inner product $\ip{\cdot}{\cdot}_X$ is as defined in \Cref{eq:ip_xv}. Then, for any $\mbf{{x}},~ \mbf{{y}}\in L_2^{n_{\mbf x}}$
	\begin{align}
	&\ip{\mcl T\mbf{x}}{\mcl{T}\mbf{y}}_X = \ip{\mbf{x}}{\mbf{y}}_{L_2}.
	\end{align}

\end{restatable}
\begin{proof}
The proof is similar to the proof of Theorem 6.3a in \cite{shivakumar_SICON2022} with $n=\{n_0,n_1,n_2\}$, $\hat{\mcl T} = \mcl T$ and other unused parameters set to empty sets or zeros.
\end{proof}

% \begin{restatable}{corollary}{corUnitary}\label{cor:unitary_T}
% Suppose $\{n,\mbf G_{\mathrm b}\}$ is admissible, $\mcl T$ is as defined in \Cref{fig:Gb_definitions} and $\mcl T$ is as defined in \Cref{fig:Gb_definitions}. If $X$ is as defined in Eq.~\eqref{eq:odepde_general_domain}, then $X$ is complete with respect to $X$-norm. Furthermore, $\mcl T: L_2 \to X$ is unitary (an isometric surjective mapping between Hilbert spaces).
% \end{restatable}
% \begin{proof}
% It is easy to show that range of an isometric (inner-product preserving) and bijective map on a complete normed-vector space is also complete and thus we proceed by just claiming without proof that $X$ is complete since the range of $\mcl T$ is $X$. Since $\mcl T$ is bijective (from \cref{thm:T_map}), preserves the inner product (from \cref{thm:unitary_T}) and is a map between Hilbert spaces, by definition, it is unitary.
% \end{proof}

% \subsection{Equivalence of norms on $\R\times X_v$ and $H$}
Now, we can show that norms induced by the inner products $\ip{\cdot}{\cdot}_{X}$ and $\ip{\cdot}{\cdot}_H$ on $X$ are equivalent and consequently, notions of stability with respect to these norms will be equivalent.

\begin{restatable}{lemma}{lemNormEquivalence}\label{lem:norm_equivalence}
Suppose pair $\{n,\mbf G_{\mathrm b}\}$ is admissible. Then, for any $\mbf{x}\in X$, $\norm{\mbf{x}}_{X}\le \norm{\mbf{x}}_H$ and there exists a constant $c_0>0$ such that $\norm{\mbf{x}}_H\le c_0\norm{\mbf{x}}_{X}$.
\end{restatable}
\begin{proof}
The proof is same as the proof of Theorem 6.4 in \cite{shivakumar_SICON2022} with $n=\{n_0,n_1,n_2\}$ and other unused parameters set to empty sets or zeros.
\end{proof}

Now that we have established that $X$-norm can be upper bounded by $H$-norm, we can use that to prove the equivalence of the stability of PDE and PIE because $X$-norm on $X$ is isometric to $L_2$-norm on $L_2$.
\begin{restatable}{theorem}{thmStability}\label{thm:stability_equivalence}
Given an $n$ and system parameters $\{\mbf G_{\mathrm b}, \mbf G_{\mathrm{p}}\}$ as defined in \Cref{eq:BC-parms,eq:pde-parms} with $\{n, \mbf G_{\mathrm b}\}$ admissible, suppose $\{\mcl{T}$ $\mcl A\}$ are as defined in \cref{fig:Gb_definitions}. Then, the PDE defined by $\{n,\mbf G_{\mathrm b}, \mbf G_{\mathrm{p}}\}$ is exponentially stable if and only if the PIE defined by $\{\mcl{T}$ $\mcl A\}$ is exponentially stable.
\end{restatable}
\begin{proof}
The proof is same as the proof of Theorem 6.6 in \cite{shivakumar_SICON2022} with $n=\{n_0,n_1,n_2\}$ and other unused parameters set to empty sets or zeros.
\end{proof}

Using the above results, we now propose the following optimization problem to test the stability of a PDE with integral terms. 

\begin{restatable}{theorem}{LPIStability}\label{thm:LPI_stability}
Given a set of PDE parameters $\{n,$ $\mbf G_{\mathrm b},$ $\mbf G_{\mathrm p}\}$, suppose there exist $\alpha, \delta >0$, $\{R_0, R_1, R_2\}$, $\{H_0,H_1,H_2\}$ such that $\threepi{R_i}, \threepi{H_i}\in [\Pi_3]_{n_{\mbf x}, n_{\mbf x}}^+$, $\threepi{R_i}\ge \alpha I$, $\threepi{H_i}\ge \delta \mcl T^*\mcl T$ and $\threepi{H_i} = -\left(\mcl T^*\threepi{R_i}\mcl A + \mcl A^*\threepi{R_i}\mcl T\right)$ where $\{\mcl{T}$, $\mcl A\}$ are as defined in \cref{fig:Gb_definitions}. Then, the PDE defined by $\{n,\mbf G_{\mathrm b}, \mbf G_{\mathrm{p}}\}$ is exponentially stable.
\end{restatable}
\begin{proof}
Suppose $R_i$ and $H_i$ are as stated above. Suppose $\mbf x$ solves the PDE defined by $\{n,$ $\mbf G_{\mathrm b},$ $\mbf G_{\mathrm p}\}$ for some initial condition $\mbf x^0\in X$. Then $\mbf x_f:= \mcl D \mbf x$ solves the PIE defined by $\{\mcl{T}$, $\mcl A\}$ for initial condition $\mcl \mbf x^0$. 

Let the Lyapunov function candidate be $V(\mbf x_f) = \ip{\mcl T\mbf x_f}{\threepi{R_i}\mcl T\mbf x_f}_{L_2}$. Then $V(\mbf x_f) \ge \alpha \norm{\mcl T\mbf x_f}_{L_2}^2$ for all $\mbf x_f\in L_2$. Taking the derivative of $V$ with respect to time along the solution trajectories of the PIE, we have
\begin{align*}
&\dot{V}(t) \\
&= \ip{\mcl T\dot{\mbf x}_f(t)}{\threepi{R_i}\mcl T\mbf x_f(t)}_{L_2}+\ip{\mcl T\mbf x_f(t)}{\threepi{R_i}\mcl T\dot{\mbf x}_f(t)}_{L_2}\\
&= \ip{\mcl A\mbf x_f(t)}{\threepi{R_i}\mcl T\mbf x_f(t)}_{L_2}+\ip{\mcl T\mbf x_f(t)}{\threepi{R_i}\mcl A\mbf x_f(t)}_{L_2}\\
&= \ip{\mbf x_f(t)}{\left(\mcl A^*\threepi{R_i}\mcl T+\mcl T^*\threepi{R_i}\mcl T\right)\mbf x_f(t)}_{L_2} \\
&= -\ip{\mbf x_f(t)}{\threepi{H_i}\mbf x_f(t)}_{L_2} \le -\delta\norm{\mcl T\mbf x_f(t)}_{L_2}^2.
\end{align*}

Then, from Gronwall-Bellman inequality, 
\[
\norm{\mbf x_f(t)}_{L_2}^2 \le \frac{k}{\alpha}\norm{\mcl D\mbf x^0}_{L_2}^2\exp(-\delta \xi t) 
\] 
where $k = \norm{\mcl T^*\threepi{R_i}\mcl T}_{\mcl L(L_2)}$ and $\xi = \norm{\mcl T}^2_{\mcl (L_2)}$. Since the initial condition was an arbitrary function, the above inequality is satisfied for any $\mcl D\mbf x^0 \in L_2$ and hence the PIE defined by $\{\mcl T, \mcl A\}$ is exponentially stable. Consequently, from \Cref{thm:stability_equivalence}, the PDE is exponentially stable.
\end{proof}

See \cite{peet_2020aut}, for a parametric form of $P_i$ and $H_i$ that allows the use of Linear Matrix Inequalities to enforce positivity constraint. Once the positivity constraint is rewritten as LMI constraints, we can use an SDP solver to find a feasible solution. 

\section{Numerical Example}\label{sec:numerical}
In this section, we present numerical tests for stability of the two PDEs introduced earlier. The steps involved in setting up the computational problem, specifically, defining the PDE parameters, converting to PIE form, setting up the operator-valued optimization problem, and converting the operator-valued optimization problem to an LMI feasibility test are all performed using the PIETOOLS toolbox for MATLAB. 
\subsection{Population Dynamics}
 Recall the McKendrick PDE model for population dynamics given by 
	\begin{align*}
	\dot{x}(t,s) &= -\partial_s x(t,s) + f(s)x(t,s),\quad s\in[0,1],\\
	x(t,0) &= \int_0^1 h(s)x(t,s)ds,
	\end{align*}
	where $x$ is the density of population of age $s$ at any given time $t$. The boundary condition can be considered as an approximation of the number of newborns at time $t$. 

	We will select the kernel in the integral term of the boundary condition as $h(s) = (1-s)s$, which implies that the population outside some normalized limits $[0,1]$ do not contribute to the birth of newborns. We will employ a constant mortality rate $f(s)=c$ and vary the $c\in \R$ to find the mortality rate, $c_0$ below which the population would go extinct (i.e, $\lim_{t\to\infty} x(t,\cdot) = 0$).

	By testing the stability of the above PDE for various $c$ values (using the method of bijection), we determined that for mortality rates greater than $-0.740625$ the population goes extinct. That is the population will survive ($x$ will not converge to zero) when the population growth rate $f(s)>0.740625$.

\subsection{Observer-based control of reaction-diffusion equation}
Consider the reaction-diffusion example presented in the introduction, where the observer is designed based on boundary measurement. We know that the observer gain, $l$, stabilizes the PDE where
\[
l(s) =-\sqrt{\lambda} \frac{I_1\left(\sqrt{\lambda(1-s^2)}\right) }{\sqrt{1-s^2}}.
\] 
Since we are interested in the practical implementation, we will approximate the gains $l$ by a polynomial of a fixed order $n$ denoted by $l_n$. Furthermore, the boundary measurements are replaced by an integral. While a typical approach is to replace point measurements by an integral with a Gaussian kernel centered at the point, in this example specifically, we can use the Fundamental Theorem of Calculus to rewrite the point measurement $\partial_s x(t,1)$ exactly as
\[\partial_s x(t,1) = \partial_s x(t,0) +\int_0^1 \partial_s^2 x(t,s) ds = \int_0^1 \partial_s^2 x(t,s) ds.\]
Likewise, we replace the point measurement value $\partial_s \hat{x}(t,1)$ by an integral to get the closed-loop observer PDE as
\begin{align*}
	&\dot{x}(t,s) = \lambda x(t,s) + \partial_s^2 x(t,s), \quad s\in [0,1],\\
	&\dot{\hat x}(t,s) = \lambda \hat x(t,s) + \partial_s^2 \hat x(t,s) \\
	&\qquad+ \int_0^1 l(s) \left(\partial_s^2 x(t,\theta) - \partial_s^2 \hat x(t,\theta)\right)d\theta,\\
	&x(t,0) = 0, \quad x(t,1) = 0,\\
	&\hat x(t,0) = 0, \quad \hat x(t,1) = 0.
	\end{align*}
Then, using the conversion formulae presented in the appendix we can find the PIE representation for the closed-loop PDE where $l$ replaced by $l_n$ which is the $n^{th}$ order polynomial approximation of $l$. 

For $\lambda\le 5$ we can prove that $1^{st}$-order polynomial approximation (a straight line) is sufficient to guarantee the stability of the closed loop PDE system. However, as $\lambda$ becomes large higher order polynomial approximation of $l$ ($l_4$ for $\lambda =6$ and so on) were necessary to prove the stability. 

\section{Conclusion}\label{sec:conclusion}
To conclude, we presented a computational method to test the stability of PDEs with integral terms appearing in the boundary conditions and dynamics where the integration is performed with respect to the spatial variable. In the process of achieving this, we presented a standard parametric form for such PDEs, provided a sufficient criterion that guarantees the existence of a PIE representation for such PDEs, and showed that both representations have the same stability properties. Finally, using the equivalence in the stability properties, we formulated the test for stability of PDEs as an optimization problem that can be solved using SDP solvers and demonstrated the application using numerical examples. 

\balance 
\bibliographystyle{plain}
\shortpage{\bibliography{references.bib}}

\begin{figure}[!h]
\centering
\scalebox{1}{
\begin{minipage}{\textwidth}
\section*{APPENDIX}
\begin{align*}
& n_{\mbf x}=\sum_{i=0}^2 n_i,\quad n_{S}=\sum_{i=1}^2 i\cdot n_{i},\quad T(s) = \bmat{I_{n_1}&(s-a)I_{n_2}\\0&I_{n_2}} \in  \R^{n_S \times n_S},\\
&U_{2i} = \bmat{0_{n_i\times n_{i+1:2}}\\ I_{n_{i+1:2}}} \in \R^{n_{S_{i}} \times n_{S_{i+1}}},\qquad  U_2 = \bmat{\text{diag}(U_{21},U_{22})\\0_{n_2\times n_S}}\in \R^{ \left(n_{\mbf x}+n_S\right)\times n_S},\\
&B_T := B\bmat{T(0)\\T(b-a)}-\int_{a}^{b}B_{I}(s)U_2 T(s-a)ds\in \R^{n_{BC}\times n_S},\quad U_{1i} = \bmat{I_{n_i} \\ 0_{n_{i+1:2}\times n_i}},\\ 
&U_1 = \bmat{U_{10}&&\\&U_{11}&\\&&U_{12}}\in \R^{(n_S+n_{\mbf x})\times n_{\mbf x}},\quad Q(s)= \bmat{0&I_{n_1}&0\\0&0&(b-s)I_{n_2}\\0&0&I_{n_2}}\in \R^{n_S \times n_\mbf x},\\
&B_Q(s) = B_T^{-1}\left(B_I(s)U_1-B\bmat{0\\Q(b-s)}+\int_s^b B_I(\theta)U_2Q(\theta-s)d\theta\right),\\
&G_0 = \bmat{I_{n_0}&\\&0_{(n_{\mbf x}-n_0)}},\qquad G_2(s,\theta) = \bmat{0\\T_1(s-a)B_Q(\theta)},\qquad G_1(s,\theta) = \bmat{0\\Q_1(s-\theta)}+G_2(s,\theta),\qquad\\
&R_{D,2}(s,\theta) = U_2T(s-a)B_Q(\theta), \qquad R_{D,1}(s,\theta) = R_{D,2}(s,\theta)+U_2Q(s-\theta), \qquad \hat{A}_0(s) = A_0(s)U_1,\\
&\hat{A}_1(s,\theta) =  A_0(s)R_{D,1}(s,\theta)+A_1(s,\theta)U_1+\int_a^{\theta} A_1(s,\beta)R_{D,2}(\beta,\theta)d\beta+\int_{\theta}^s A_1(s,\beta)R_{D,1}(\beta,\theta)d\beta \\&\quad + \int_{s}^b A_2(s,\beta)R_{D,1}(\beta,\theta)d\beta,\\
&\hat{A}_2(s,\theta) =  A_0(s)R_{D,2}(s,\theta)+A_2(s,\theta)U_1+\int_a^{s} A_1(s,\beta)R_{D,2}(\beta,\theta)d\beta+\int_s^{\theta} A_2(s,\beta)R_{D,2}(\beta,\theta)d\beta \\&\quad + \int_{\theta}^b A_2(s,\beta)R_{D,1}(\beta,\theta)d\beta,\\
& \mcl{T} = \threepi{G_i},\qquad \mcl A = \threepi{\hat A_i}.
\end{align*}
\caption{Definitions based on $\{n,A_i,B,B_I\}$.}\label[EqnBlock]{fig:Gb_definitions}
\end{minipage}}

\end{figure}

	\end{document}